\documentclass[12pt,twoside]{article}
\usepackage{amsmath}
\usepackage{amsfonts}
\usepackage{amsthm}
\usepackage{amssymb}
\usepackage{graphicx}
\usepackage{multicol}

\begin{document}
\title{{\normalsize{\bf UNIQUE DIAGRAM OF A SPATIAL ARC AND THE KNOTTING 
PROBABILITY}}}
\author{{\footnotesize Akio KAWAUCHI}\\
{\footnotesize{\it Osaka City University Advanced Mathematical Institute}}\\
{\footnotesize{\it Sugimoto, Sumiyoshi-ku, Osaka 558-8585, Japan}}\\
{\footnotesize{\it kawauchi@sci.osaka-cu.ac.jp}}}
\date\, 
\maketitle
\vspace{0.25in}
\baselineskip=10pt
\newtheorem{Theorem}{Theorem}[section]
\newtheorem{Conjecture}[Theorem]{Conjecture}
\newtheorem{Lemma}[Theorem]{Lemma}
\newtheorem{Sublemma}[Theorem]{Sublemma}
\newtheorem{Proposition}[Theorem]{Proposition}
\newtheorem{Corollary}[Theorem]{Corollary}
\newtheorem{Claim}[Theorem]{Claim}
\newtheorem{Definition}[Theorem]{Definition}
\newtheorem{Example}[Theorem]{Example}

\begin{abstract} It is shown that the projection image of an oriented spatial arc to any oriented plane is approximated  by a 
unique arc diagram (up to isomorphic arc diagrams) determined from the spatial arc 
and the projection. 
In a separated paper, the knotting probability of an arc 
diagram is defined as  an invariant under isomorphic arc diagrams. By combining them, the knotting probability of every oriented spatial arc is defined.  

\phantom{x}


\end{abstract}

\baselineskip=15pt

\bigskip

\noindent{\bf 1. Introduction}

A  {\it spatial arc} is a  polygonal arc in the 3-space  ${\mathbf R}^3$, which is  considered as a model of a protein or a linear polymer in science. 
The following question on science is an interesting question that can be set as a mathematical question:

\phantom{x}

\noindent{\bf Question.} 
How a linear scientific object such as a linear molecule (e.g. a non-circular DNA, protein, linear polymer,…) is considered as a knot object ?

\phantom{x}

In this paper,  it is shown that the projection image of an oriented spatial arc to any oriented plane is approximated into a 
unique arc diagrams (up to isomorphic arc diagrams) determined from  
the spatial arc and the projection. 
Further, the orientation change of the  spatial arc 
makes only substitutes for these unique arc diagrams (up to isomorphic arc diagrams). 
This argument is more or less similar to an argument transforming a classical knot 
in ${\mathbf R}^3$ into a regular knot diagram (see \cite{CF}, \cite{K0}). 
Let
\[S^2=\{ u\in  {\mathbf R}^3 |\, ||u||=1\}\]
be the unit sphere, where $||\phantom{u}||$ denotes the norm on ${\mathbf R}^3$. 
Every element $u\in S^2$ is regarded as a unit vector from the origin $0$.
For a unit vector $u\in {\mathbf R}^3$, let $P_u$ be the oriented plane containing 
the origin $0$ such that the unit vector $u$ is a positively normal vector to $P_u$. 
The orthogonal projection from ${\mathbf R}^3$ to the plane $P_u$ is called the 
{\it projection along} the unit vector $u\in S^2$ and denoted by 
\[\lambda_u:{\mathbf R}^3\to P_u.\] 

For a small positive number $\delta$, a $\delta$-{\it approximation} of 
the projection $\lambda_u:{\mathbf R}^3\to   P_u$ along $u\in S^2$  is the projection 
\[\lambda_{u'}:{\mathbf R}^3\to   P_{u'}\] 
along a unit vector $u'\in S^2$  with $||u'-u||<\delta$, 
which is denoted by 
\[\lambda_{u}^{\delta}: {\mathbf R}^3\to P_u.\]
The projection image $\lambda_u(L)$ of  a spatial arc $L$  in  the plane $P_u$
is an {\it  arc diagram} in the oriented plane $P_u$ 
if $\lambda_u(L)$ has only {\it crossing points} 
(i.e., transversely meeting double points with over-under information) and 
the starting and terminal points  as single points. 
The arc diagram $(P_u, \lambda_u(L))$ is sometimes denoted by $(P,D)$.

An arc diagram $(P,D)$  is {\it isomorphic} to an arc diagram $(P',D')$  if there 
is an orientation-preserving homeomorphism $f:P\to P'$ sending $D$ to $D'$ 
which  preserves the crossing points of $D$ and $D'$, and the starting and terminal points of $D$ and $D'$. The map $f$ is called an 
{\it isomorphism} from $D$ to $D'$. 
In an illustration of an arc diagram, it is convenient to illustrate 
an arc diagram with smooth edges 
in the class of isomorphic arc diagrams instead of a polygonal arc diagram. 

In this paper, the following observation is shown.

\phantom{x}

\noindent{\bf Theorem~1.1.}  
Let $L$ be an oriented arc  in ${\mathbf R}^3$, and  
$\lambda_u:{\mathbf R}^3\to P_u$  the projection 
along a unit vector $u\in S^2$. 
For any sufficiently small positive number $\delta$,  
the projection  $\lambda_u$ has a $\delta$-approximation
\[\lambda_{u}^{\delta}: {\mathbf R}^3\to P_{u}\] 
such that the projection image $\lambda_u^{\delta}(L)$ 
is an arc diagram  determined uniquely 
from  the spatial arc $L$ and the projection $\lambda_u$ 
up to isomorphic arc diagrams. 

\phantom{x}

The proof of Theorem~1.1 is done in \S~2.
In \cite{K7}, the knotting probability 
\[p(D)=(p^{\tiny\mbox{I}}(D),p^{\tiny\mbox{II}}(D),p^{\tiny\mbox{III}}(D))\] 
of an arc diagram $D$ is defined so that 
it is unique up to isomorphic arc diagrams.  
By the arc diagram $D(L;u)=\lambda_u^{\delta}(L)$,  
the {\it knotting probability} $p(L;u)$ of an oriented spatial arc $L$ is 
defined by 
\[p(L;u)=p(D(L;u))\] 
for every unit vector $u\in S^2$.  More details are discussed in \S~3. 

We mention here that a knotting probability of a spatial arc was defined directly from a knotting structure of a spatial graph but with the demerit that it depends on 
the heights of the crossing points of a diagram of the spatial arc in \cite{K1, K3}.

\phantom{x}

\noindent{\bf 2. Proof of Theorem~1.1}

Let $L$ be an oriented spatial arc, and  
$p_s$ and $p_t$ the starting point  and the terminal point of $L$, respectively. 
The {\it front edge} of $L$ is the interval $\gamma$ in ${\mathbf R}^3$ joining    
the starting point $p_s$ and the terminal point $p_t$. 
Orient $\gamma$ by the orientation from $p_s$ to $p_t$. 
Let $u({\gamma})$ be the unit vector of the front edge $\gamma$ of $L$ called the 
{\it front vector} of $L$. 

Let $L$ be an oriented spatial arc with  the starting point  
$p_s$ and the terminal point $p_t$. 
An {\it edge line} of $L$ is an oriented line\footnote{Throughout the paper, by a 
{\it line},  
we mean a straight line  $\ell$ in ${\mathbf R}^3$ 
extending an edge of $L$ oriented by $L$.}  
The {\it front edge} of $L$ is the interval $\gamma$ in ${\mathbf R}^3$  joining    
the starting point $p_s$ and the terminal point $p_t$. 
Orient $\gamma$ by the orientation from $p_s$ to $p_t$. 
The {\it front line} $\ell_{\gamma}$ is the oriented line extending 
the oriented front edge $\gamma$ of $L$  from     
the starting point $p_s$ to the terminal point $p_t$. 

Assume that there is an edge line of $L$ distinct from the front line $\ell_{\gamma}$ 
because otherwise  there is nothing to show.   
The {\it starting front-pop line} of $L$ is the edge line  $\ell_s$  of  the edge 
which pops for the first time from the front line $\ell_{\gamma}$ 
when a point is going on $L$ along  the orientation of $L$. 
The {\it ending front-pop line} of $L$ is the edge line $\ell_t$ of the edge which reaches the front line $\ell_{\gamma}$  at the end 
when a point is going on $L$ along the orientation of $L$. 

The {\it starting front-pop plane} of an oriented spatial arc $L$ is the oriented plane $P(\ell_{\gamma},\ell_s)$ determined by the front line $\ell_{\gamma}$ and the starting pop line $\ell_s$ in this order. 
The {\it terminal front-pop plane} of an oriented spatial arc $L$ is the oriented plane $P(\ell_{\gamma},\ell_t)$ determined by the front line $\ell_{\gamma}$ and the terminal  pop line $\ell_s$ in this order. 

Let $u(\ell)\in S^2$ be the unit vector of  an oriented line $\ell$ in ${\mathbf R}^3$. 
Let $u_{\gamma}\in S^2$ denote the unit vector of $\ell_{\gamma}$, called 
the {\it front edge vector}. 
Let $u_s, u_t\in S^2$ denote the unit vectors of the starting and terminal front-pop line $\ell_s$ and $\ell_t$, called the 
{\it starting and terminal front-pop}  vectors, respectively. 

For a plane $P$ in ${\mathbf R}^3$, the {\it great circle} $C$  of $P$ in $S^2$
is  the great circle obtained as the intersection of $S^2$ and a plane $P'$ 
parallel to $P$.  

The {\it trace set} $T$ of a spatial arc $L$ is the subset of $S^2$ consisting  
of  the great circles and the unit vectors obtained from  $L$  in the following 
cases  (i) and (ii)* 

\phantom{x}

\noindent{(i)} The great circle $C$ of  $S^2$  of the plane $P$ 
in ${\mathbf R}^3$ determined by a vertex  $u$ of  $L$  
and an edge line $\ell$ or the front line $\ell_{\gamma}$ of $L$ which is 
disjoint from $u$. 

\medskip

\noindent{(ii)} The  unit vectors $\pm u_{\eta}\in S^2$  of  an oriented line 
$\eta$ in ${\mathbf R}^3$ meeting three edge lines  $\ell_i\, (i=1,2,3)$  of $L$  
any two of which are not on the same plane with $3$ distinct points.

\phantom{x}

By (i), note that the following great circles are in the trace set $T$: 

\phantom{x}

\noindent(1) the great circles of the starting and terminal front-pop planes  $P(\ell_{\gamma},\ell_s)$ and $P(\ell_{\gamma},\ell_t)$, 

\medskip

\noindent(2) the great circle of the plane determined by  two parallel distinct edge lines $\ell and\ell'$, 

\medskip

\noindent(3) the great circle of the plane determined by two distinct edge lines 
$\ell, \ell'$ meeting a point, 

\medskip

\noindent(4)  the great circle of the plane determined by an edge line $\ell$ 
and the front line $\ell_{\gamma}$ meeting a point.   

\phantom{x}

In particular, the unit vectors $\pm u(\ell)\in S^2$ of every edge line $\ell$ of $L$  
and the unit vectors $\pm u_{\gamma}\in S^2$  and 
the great circle of any plane  containing any three distinct lines 
$\ell_i\, (i=1,2,3)$ is the trace set $T$. 
Also the  unit vectors $\pm u(\eta)\in S^2$  of  an oriented line 
$\eta$ in ${\mathbf R}^3$ meeting three edge lines  $\ell_i\, (i=1,2,3)$  of $L$  
some two of which are on the same plane with $3$ distinct points.  

In (ii), note that a line $\eta$ meeting $\ell_i\, (i=1,2,3)$  by different 3 points 
is unique, because 
if there is another such line $\eta'$, then the lines  
$\eta, \eta', \ell_i\,(i=1,2,3)$ 
and hence the lines  $\ell_i\, (i=1,2,3)$ are on the same plane, contradicting the assumption.  

Also, note that  if a unit vector $u\in S^2$ is in the trace set $T$, then 
$-u$ is also in $T$.

For every unit vector $u\in S^2$, let $P_u$ be the oriented plane containing 
the origin $0$ such that $u$ is positively normal to $P_u$, 
and $\lambda_u:{\mathbf R}^3 \to P_u$ the orthogonal projection.  
We show the following lemma. 

\phantom{x}

\noindent{\bf Lemma~2.1.} For every unit vector $u\in S^2\setminus T$, 
the projection image $\lambda_u(L)$ is an arc diagram in the plane $P_u$.  
Further, the arc diagram  $\lambda_u(L)$ up to isomorphic arc diagrams is independent 
of  any choice of a unit vector $u'$ in the connected region $R(u)$ of $S^2\setminus T$ 
containing $u$. 

\phantom{x}

\noindent{\bf Proof of Lemma~2.1.}
If a unit vector $u\in S^2$ is not in (i), then every edge of $L$ and  the front line 
$\ell_{\gamma}$ are embedded in the plane $P_u$ by the projection $\lambda_u$. 
If a unit vector $u\in S^2$ is  in neither (i) nor (ii), then the set of vertexes 
of $L$ is embedded in $P_u$ by the projection $\lambda_u $   whose image is disjoint 
from the image of any open edge of $L$. 
In particular, any two distinct parallel edge lines are disjointedly 
embedded in $P_u$. 
Further, the images of the edges of $L$ meet only in the images of the open edges 
of $L$. 
Thus, if a unit vector $u\in S^2$ is in neither (i) nor (ii), namely if 
$ u\in S^2\setminus T$ , then the meeting points among the edges of $L$ 
consisting of double points between two open edges of $L$ and hence 
the projection image $\lambda_u(L)$ is an arc diagram in the plane $P_u$. 
The arc diagram  $\lambda_{u'}(L)$ is unchanged up to isomorphisms 
for any unit vector $u'$ in a connected open neighborhood of $u$ in $S^2\setminus T$, 
so that the arc diagram  $\lambda_{u'}(L)$ is unchanged up to isomorphisms 
for any unit vector $u'$ in the connected region $R(u)$. 
$\square$

\phantom{x}

The proof of Theorem~1.1 is done as follows:

\phantom{x}

\noindent{\bf Proof of Theorem~1.1.}
The idea of the proof is to specify a unique connected region 
$R(u')$ of $S^2\setminus T$ adjacent to every unit vector $u\in S^2$. 
For this purpose,  for any given oriented spatial arc $L$ (not in the front line 
$\ell_{\gamma}$), the new $x$-axis, $y$-axis and $z$-axis of the 3-space 
${\mathbf R}^3$ are set as follows:

The unit vector $u_{\gamma}$ of the front edge $\gamma$ is taken as the unit vector of the $x$-axis as follows: 
\[e_x=u_{\gamma}=(1,0,0).\]
Let $y_s$ be the unit vector modified from the starting front-pop vector $u_s$ 
to be  orthogonal to the unit vector 
$e_x=u_{\gamma}$   with the inner product $y_s\cdot u_s>0$ 
in the starting front-pop plane $P(\ell_{\gamma},\ell_s)$ of $L$. 
The unit vector $y_s$ is taken as the unit vector the $y$-axis:
\[e_y=y_s=(0,1,0).\]
Then exterior product $z_s= u_{\gamma}\times y_s$ which is given by the z-axis:
\[e_z=z_s=(0,0,1).\]

Under this setting of the coordinate axis, let $S^2$ be the unit sphere 
which is the union of the upper hemisphere $S^2_+$, the equatorial circle $S^2_0$ 
and the lower hemisphere $S^2_-$ given as follows:
\begin{eqnarray*}
S^2_+&=&\{(x,y,z)\in {\mathbf R}^3|\, x^2+y^2+z^2=1, \, z>0\},\\
S^2_0&=&\{(x,y,z)\in {\mathbf R}^3|\, x^2+y^2+z^2=1, \, z=0\},\\
S^2_-&=&\{(x,y,z)\in {\mathbf R}^3|\, x^2+y^2+z^2=1, \, z>0\}.
\end{eqnarray*}
Note that the equatorial circle $S^2_0$ belongs to $T$  since the front-pop plane $P_L$  
coincides with the plane with $z=0$. 

Every unit vector $u$ in $S^2_+$ except the north pole 
$(0,0,1)$ is uniquely written as 
\[u=\psi(r,\theta)=(r\cos \theta, r\sin \theta, \sqrt{1-r^2})\] 
for real numbers $r$ and $\theta$ with $0<r\leq 1$ and $0\leq \theta<2\pi$  
in a unique  way. 

\phantom{x}  

\noindent{Case 1}: $u=\psi(r,\theta)\in S^2_+$ with $0<r<1$. 

Note that when the number $r$ with $0<r<1$ is fixed,  
the points $\psi(r,\theta')$ for all $\theta'$ with  $0\leq \theta'<2\pi$  form 
a circle in $S^2$ which is different from every great circle of $S^2$ 
and hence meets $T$  only in finitely many points.  

The connected region $R(u')$ is taken   
to be  the connected region  of $S^2_+\setminus T$ which is adjacent 
to the unit vector $\psi(r,\theta')$ and contains the 
unit vector $\psi(r, \theta+\varepsilon)$ for a sufficiently small 
positive number $\varepsilon$, which is uniquely determined.

\medskip

\noindent{Case 2}: $u=(0,0,1)\in S^2_+$. 
By taking a positive number $r$ sufficiently small, 
take a unit vector $\psi(r, 0)=(r,0,\sqrt{1-r^2})\in S^2_+$  
such that $\psi(r', 0)$ does not meet the great circles in $T$ for any $r'$ with 
$0<r'\leq r$  except for the meridian  circle $x^2+z^2=1$ (if it is in $T$). 
Then the connected region $R(u')$ is taken to be 
the connected region of $S^2_+\setminus T$ which is adjacent to the unit vector 
$\psi(r, 0)$ and contains the unit vector 
$\psi(r, \varepsilon)$ for a sufficiently small positive number 
$\varepsilon$, which is uniquely determined.

\medskip

\noindent{Case 3}: $u=\psi(1,\theta)\in S^2_0$.

Let $0\leq \theta<\pi$. By taking a positive number $r$ with $1-r$ a sufficiently small positive number, take a unit vector $\psi(r,\theta)\in S^2_+$ 
such that $\psi(r',\theta)$  does not meet the great circles in $T$ for any $r'$ 
with $r\leq r'<1$.  
Then the connected region $R(u')$ is taken to be the connected region 
of $S^2_+\setminus T$ which is adjacent to the unit vector 
$\psi(r, \theta)$ and contains the unit vector 
$\psi(r, \theta+\varepsilon)$ for a sufficiently small 
positive number $\varepsilon$, which is uniquely determined. 

Let $\pi\leq \theta<2\pi$. Since $-u=\psi(1,\theta-\pi)\in S^2_0$,  
we specified the connected component 
$R((-u)')$ of $S^2\setminus T$. The  desired connected component 
$R(u')$ of $S^2\setminus T$  is the image  of $R((-u)')$ under the antipodal map 
$-1:S^2\to S^2$ defined by  $(x,y,z)$ to $(-x,-y,-z)$ with $(-1)(T)=T$. 

\medskip

\noindent{Case 4}: $u\in S^2_-$.  

Since the unit vector $-u$ is in $S^2_+$, we specified the connected component 
$R((-u)')$ of $S^2\setminus T$ in the cases 1-3. The  desired connected component 
$R(u')$ of $S^2\setminus T$  is the image  of $R(u')$ under the antipodal map 
$-1:S^2\to S^2$.  

\phantom{x}

Thus,  for every unit vector $u\in S^2$, the connected region $R(u')$ and  
is  specified, so that for every $u''\in R(u')$, 
the image $\lambda_{u''}(L)$ of $L$ by the orthogonal projection 
$\lambda_{u''}:{\mathbf R}^3 \to P_{u''}$ 
is an arc diagram determined uniquely up to isomorphic arc diagrams from $L$ and 
the projection $\lambda_u$.  

Since the connected region $R(u')$ is  adjacent to the normal vector $u\in S^2$,  
for every $\delta>0$  there is a unit vectors $u''\in R(u')$   
with $||u''-u'||<\delta$ and  by Lemma~3.1, the orthogonal projection 
$\lambda_{u''}:{\mathbf R}^3 \to P_{u''}$ is a desired $\delta$-approximation   
$\lambda_u^{\delta}: {\mathbf R}^3\to P_u$. 
This completes the proof of Theorem~1.1. $\square$

\phantom{x}  

An arc diagram is {\it inbound} if the starting and terminal points are in the same region of the plane divided by the arc diagram. 
The  projection image $\lambda_u(L)$ is {\it inbound} if 
the image of the front edge $\gamma$ of $L$ does not meet 
the image of any edge of $L$ under the projection $\lambda_u$ 
except for the starting and terminal points. 
A spatial arc $L$ is {\it inbound} if the interior of the front edge $\gamma$ of  
an oriented spatial arc $L$ does not meet $L$. 
A spatial arc $L$ is {\it even} if the starting front-pop plane 
$P(\ell_{\gamma},\ell_s)$ and the terminal front-pop plane $P(\ell_{\gamma},\ell_t)$ 
coincide as a set. 
Let $-L$ denote the same spatial arc as $L$ but with the opposite orientation. 
We have the following observations from the proof of Theorem~1.1. 

\phantom{x}

\noindent{\bf Corollary~2.2.} For an oriented spatial arc $L$ and an arc 
diagram 
 $D(L;u)=\lambda_{u}^{\delta}(L)$,   we have the following (1)-(5).

\phantom{x}

\noindent{(1)} The arc diagram $D(L;-u)$ is the  mirror image of $D(L;u)$.  

\medskip

\noindent{(2)} If the projection image $\lambda_u(L)$ is an arc diagram, then   
the arc diagram $D(L;u)$ is isomorphic to the  arc diagram $\lambda_u(L)$.

\medskip

\noindent{(3)}  There are  only finitely many arc diagrams $D(L;u)$ 
up to isomorphic arc diagrams for  all unit vectors $u\in S^2$.  

\medskip

\noindent{(4)} If the projection image  $\lambda_u(L)$ is inbound, then   
the arc diagram $D(L;u)$ is an inbound arc diagram.
In particular,  
if the spatial arc $L$ is inbound, then 
the arc diagram $D(L;u_{\gamma})$ for the front edge vector 
$u_{\gamma}$ is an inbound arc diagram.

\medskip

\noindent{(5)}  If the spatial arc $L$ is even, then 
the arc diagram $D(-L;u)$ is isomorphic to the arc diagram $D(L;u)$ 
or the mirror image $D(L;-u)$ of $D(L;u)$ with the string orientation changed 
according to whether the orientations of $P(\ell_{\gamma},\ell_s)$ and  $P(\ell_{\gamma},\ell_t)$ coincide or not. 

\phantom{x}

\noindent{\bf3. The knotting probability of a spatial arc}

A {\it chord graph} is a trivalent connected graph $(o;\alpha)$ in 
 ${\mathbf R}^3$ consisting of a trivial oriented link $o$ (called a 
{\it based loop system}) and the attaching arcs $\alpha$ 
(called a {it chord system}),
where some chords of $\alpha$  may meet.
A {\it chord diagram} is a diagram $C(o;\alpha)$ (in a plane) of a 
chord graph $(o;\alpha)$.  

A {\it ribbon surface-link}  is a surface-link in the 4-space ${\mathbf R}^4$ 
obtained from a trivial $S^2$-link by mutually disjoint embedded 1-handles (see 
\cite[II]{KSS}, \cite{Yanagawa}). 
From a chord diagram $C$, a ribbon surface-link $F(C)$ in the 4-space 
${\mathbf R}^4$ is constructed so that 
an equivalence of $F(C)$ corresponds to a combination of  
the moves $M_0$, $M_1$ and $M_2$  (see \cite{ K2,K4,K5,K6}). 
 
\begin{figure}[hbt]
\centerline{\includegraphics[width=7.5 cm,height=2.5 cm]{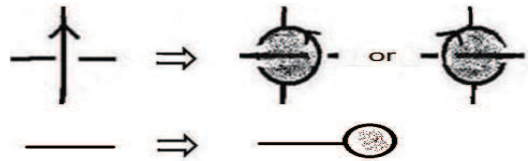}}
\caption{Transformation of an oriented arc diagram into a chord diagram}
\label{fig:trans}
\end{figure}

Let $D$  be an oriented $n$-crossing arc diagram, and $C(D)$ 
the chord diagram of $D$.  
We obtain a chord diagram $C(D)$ 
with $n+2$ based loops  from the diagram $D$ by replacing every crossing 
point and every endpoint with a based loop as in Fig.~\ref{fig:trans}.

\begin{figure}[hbt]
\centerline{\includegraphics[width=12 cm,height=4.5 cm]{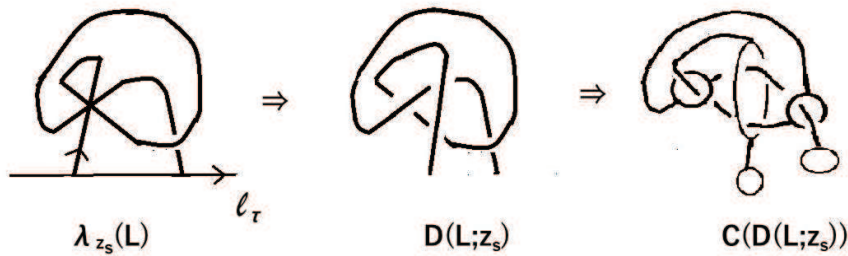}}
\caption{Producing the chord diagram $C(D(L;z_s))$ from 
the projection image $\lambda_{z_s}(L)$ }
\label{fig:DIrreg}
\end{figure}

If two based loops are connected by a chord not meeting the other chords, we can replace it by one based loop (which 
is also called the {\it chord diagram} of an arc diagram)  although 
a carefulness is needed in the calculation of the knotting probability. 
An example of the procedure to obtain the chord diagram $C(D(L;z_s))$ from 
the projection image $\lambda_{z_s}(L)$ of a spatial arc $L$ with  $3$ edges of the heights suitable near the triple point  is illustrated 
in Fig.~\ref{fig:DIrreg}, where note that the case of $u=z_s$ is easier to construct the chord diagram $C(D(L;u))$ from the projection image $\lambda_u(L)$ of a spatial arc $L$. 

Let $D$  be an oriented $n$-crossing arc diagram, and $C(D)$ 
the chord diagram of $D$.  Let $o_s$ and $o_t$ be the based loops in $C(D)$ 
transformed from the starting and terminal points $v_s$ and $v_t$, respectively.
There are $(n+2)^2$ chord diagrams $A$ obtained from the chord diagram $C(D)$ 
by joining the loops $o_s$ and  $o_t$  with any based loops of $C(D)$ by two chords 
not passing the other based loops.  
A chord diagram obtained in this way is called an  {adjoint chord diagram} of $C(D)$ 
with an {\it additional chord pair}. 
Note that the ribbon surface-knot $F(C(D))$ of the chord diagram $C(D)$  
is a ribbon $S^2$-knot and the ribbon surface-knot $F(A)$ of 
an adjoint chord diagram $A$  is a genus $2$ ribbon surface-knot. 
A chord diagram is said to be {\it unknotted} or {\it knotted} according to whether  
it represents a trivial or non-trivial ribbon surface-knot, respectively. 

The idea of the knotting probability is to measure how many knotted chord diagrams there are among the $(n+2)^2$ adjoint chord diagrams  of $C(D)$. 
Since there are overlaps among them up to  canonical isomorphisms, we consider the $n^2+2n+2$ adjoint chord diagrams $A$ of $C(D)$ by removing some overlaps. 
The $n^2+2n+2$ adjoint chord diagrams $A$ of $C(D)$ are classified by the following three types:

\phantom{x}

\noindent{\bf Type I.} Here are the $2$ adjoint chord diagrams of $C(D)$ which 
are  the adjoint chord diagram 
with two self-attaching additional chords and 
the adjoint chord diagram with a self-attaching additional 
chord on $o_s$ and an additional chord joining $o_s$ with $o_t$. 

\medskip 

\noindent{\bf Type II.} Here are the $2n$ adjoint chord diagrams $A$ of $C(D)$. 
The $2n$ adjoint chord diagrams of $C(D)$ are given by 
the additional chord pairs consist of a self-attaching additional chord on 
$o_s$ (or $o_t$, respectively) and an additional chord joining 
$o_t$ (or $o_s$, respectively) with a based loop except for $o_s$ and $o_t$. 

\medskip 

\noindent{\bf Type III.} Here are the $n$ adjoint chord diagrams $A$ of $C(D)$ where 
the additional chord pairs consist of an additional chord joining  $o_s$ with  
$o_t$  and an  additional chord joining $o_s$ with a  
a based loop except for  $o_s$  and  $o_t$.

\medskip 

\noindent{\bf Type IV.} Here are the $n(n-1)$ adjoint chord diagrams $A$ of $C(D)$ where 
the additional chord pair joins the pair of $o_s$ and $o_t$ with a distinct based loop pair not containing $o_s$ and $o_t$. 
\phantom{x}

\phantom{x}

In \cite{K7}, it is shown that every adjoint chord diagram of the chord diagram 
$C(D)$ of any $n$ crossing arc diagram $D$  is deformed into one of 
the adjoint chord diagrams of type I, II, III and IV of the chord diagram $C(D)$.
Thus,  it is justified to reduce   
the $(n+2)^2$ adjoint chord diagrams to the $n^2+2n+2$ adjoint chord diagrams.

The {\it knotting probability} $p(D)$ of an arc diagram $D$ is defined to be the quadruple 
\[p(D)=(p^{\tiny\mbox{I}}(D),p^{\tiny\mbox{II}}(D), 
p^{\tiny\mbox{III}}(D), p^{\tiny\mbox{IV}}(D))\]
of  the following knotting probabilities $p^{\tiny\mbox{I}}$, 
$p^{\tiny\mbox{II}}$, $p^{\tiny\mbox{III}}(D)$ of types I, II III, IV.

\phantom{x}

\noindent{\bf Definition.} 

\noindent{(1)} Let  $A_1$ and $A_2$ be the adjoint chord diagrams of type I  and assume that there are just $k$ knotted chord diagrams among them. Then 
the {\it type I knotting probability} of $D$  is        
\[p^{\tiny\mbox{I}}(D)=\frac{k}{2}.\]

\medskip 

\noindent{(2)} Let  $A_i\,(i=1,2,\dots, 2n)$ be the adjoint chord diagrams of 
type II  and 
assume that there are just $k$ knotted chord diagrams  among them. Then 
the {\it type II knotting probability} of $D$  is        
\[p^{\tiny\mbox{III}}(D)=\frac{k}{2n}.\]

\medskip 

\noindent{(3)} Let  $A_i\,(i=1,2,\dots, n)$ be the adjoint chord diagrams of type III  and assume that there are just $k$ knotted chord diagrams  among them. Then 
the {\it type III knotting probability} of $D$  is        
\[p^{\tiny\mbox{III}}(D)=\frac{k}{n}.\]

\medskip 

\noindent{(4)} Let  $A_i\,(i=1,2,\dots, n(n-1))$ be the adjoint chord diagrams of 
type IV  and 
assume that there are just $k$ knotted chord diagrams  among them. Then 
the {\it type IIV knotting probability} of $D$  is        
\[p^{\tiny\mbox{IV}}(D)=\frac{k}{n(n-1)}.\]

\phantom{x} 

When the orientation of an arc diagram $D$  is changed,  all the orientations of the based loops  of the chord graph $C(D)$ are changed at once. 
This means that the knotting probability $p(D)$  does not depend on any choice of orientations of $D$, and we can omit the orientation of $D$ in figures. 
See \cite{K7} for actual calculations of $p(D)$.

The knotting probability $p(D)$ has $p(D)=1$ if 
\[p^{\tiny\mbox{I}}(D)=p^{\tiny\mbox{II}}(D)
=p^{\tiny\mbox{III}}(D)=p^{\tiny\mbox{IV}}(D)=1\]
and otherwise, $p(D)<1$. 
The knotting probability $p(D)$ has $p(D)>0$ if 
\[p^{\tiny\mbox{I}}(D)+p^{\tiny\mbox{II}}(D)
+p^{\tiny\mbox{III}}(D)+p^{\tiny\mbox{IV}}(D)>0\]
and otherwise, $p(D)=0$.

For an oriented spatial arc $L$ in ${\mathbf R}^3$,   
the {\it knotting probability} $p(L;u)$ of  $L$ for a unit vector 
$u\in S^2$ is defined to be 
\[p(L;u)= p(D(L;u))\]
for the arc diagram $D(L;u)$ for $u$. Thus, 
\begin{eqnarray*}
p^{\tiny\mbox{I}}(L;u) &=& p^{\tiny\mbox{I}}(D(L;u)),\\
p^{\tiny\mbox{II}}(L;u) &=& p^{\tiny\mbox{II}}(D(L;u)),\\ 
p^{\tiny\mbox{III}}(L;u) &=& p^{\tiny\mbox{III}}(D(L;u)),\\
p^{\tiny\mbox{IV}}(L;u) &=& p^{\tiny\mbox{IV}}(D(L;u)).
\end{eqnarray*}

Although $p(L;u)$ is unchanged by any choice of a string orientation in the 
arc diagram level $D(L;u)$, the knotting probability $p(-L;u)$ may 
be much different from the arc diagram $p(L;u)$ in general,  
because the arc diagram $D(L;u)$ may be much different from the arc diagram $D(-L;u)$ 
in general (cf. Corollary~2.4).
In any case, the unordered pair of $p(L;u)$ and $p(-L;u)$ is considered as the knotting probability of a spatial arc $L$ independent of the string orientation. 
When one-valued probability is desirable, 
a suitable average of the knotting probabilities 
\[p^{\tiny\mbox{I}}(\pm L;u),\, p^{\tiny\mbox{II}}(\pm L;u),\, 
p^{\tiny\mbox{III}}(\pm L;u),\, 
p^{\tiny\mbox{IV}}(\pm L;u)\] 
is considered. 
For a special spatial arc $L$, we have the following corollary. 

\phantom{x}

\noindent{\bf Corollary~3.1.} 

\noindent{(1)} If the projection image  $\lambda_u(L)$ is inbound, then   
\[p(L;u)=p(L;-u)\] 
for any unit vector $u\in S^2$. In particular, 
if the spatial arc $L$ is inbound, then 
$p(L;u_{\gamma})=p(L;-u_{\gamma})$ for the front edge vector $u_{\gamma}$.

\medskip

\noindent{(2)}  If  $L$ is even, then 
\[p(-L;u)=p(L;u)\quad \mbox{or}\quad p(-L;u)=p(L;-u)\] 
according to whether the orientations of $P(\ell_{\gamma},\ell_s)$ and  $P(\ell_{\gamma},\ell_t)$ coincide or not. 

\phantom{x}

\noindent{\bf Proof of Corollary~3.1.} In \cite{K7},  it is shown 
$p(D)=p(D^*)$ for an inbound arc diagram $D$ and the mirror image $D^*$ of $D$.  
By Corollary~2.2 (1) and (4), the arc diagrams $D(L;u)$ and $D(L;u_{\gamma})$ are 
inbound arc diagram with $D(L;-u)$ and $D(L;-u_{\gamma})$ the mirror images. 
Thus, (1) is obtained. For (2), by Corollary~2.2 (1) and (5),
the arc diagram $D(-L;u)$ is isomorphic to the arc diagram $D(L;u)$ or $D(L;-u)$ according to whether the orientations of $P(\ell_{\gamma},\ell_s)$ and  $P(\ell_{\gamma},\ell_t)$ coincide or not, obtaining the result. 
$\square$

\phantom{x}

\begin{figure}[hbt]
\centerline{\includegraphics[width=12 cm,height=8cm]{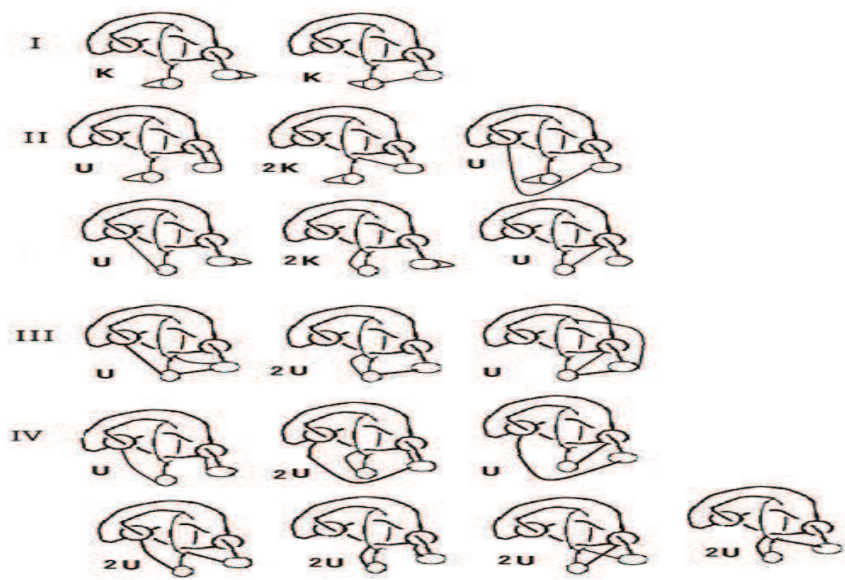}}
\caption{Calculations on the chord diagram $C(D(L;z_s))$ }
\label{fig:DIrregCal}
\end{figure}

For example, the spatial arc $L$  in Fig.~\ref{fig:DIrreg} is inbound and even, and
we have 
\[p(L;z_s)=p(L;-z_s)=p(-L,z_s) =p(-L,-z_s)=(1, \frac{1}{2},0,0)\]
by Corollary~3.1 where the calculations on the arc diagram $C(D(L;z_s))$  
are done in Fig.~\ref{fig:DIrregCal} (see \cite{K7} for the details of the calculation).

\phantom{x}

\noindent{\bf Acknowledgements.} This work was in part supported by Osaka City University Advanced Mathematical Institute (MEXT Joint Usage/Research Center on  Mathematics and Theoretical Physics).

\phantom{x}

\end{document}